\documentstyle[11pt]{article}

\newtheorem{theorem}{{\bf Theorem}}[section]

\newtheorem{lemma}{{\bf Lemma}}[section]

\newtheorem{proposition}{{\bf Proposition}}[section]

\newtheorem{definition}{{\bf Definition}}[section]

\newtheorem{corollary}{{\bf Corollary}}[section]

\newtheorem{remark}{{\bf Remark}}[section]

\def\O{\Omega}  

\def\P{{\bf P}}

\def\E{{\bf E}}

\def\F{{\cal F}}

\def\N{{\rm I\!N}}

\def\H2{{\rm I\!H2}}

\def\R{{\rm I\!R}}

\def\toD{{\scriptstyle\,\buildrel{\D}\over{\hbox to
0.6cm{\rightarrowfill}}\,}}

\def\toP{{\scriptstyle\,\buildrel{\P}\over{\hbox to
0.6cm{\rightarrowfill}}\,}}

\def\toL{{\scriptstyle\,\buildrel{\cal L}\over{\hbox to
0.6cm{\rightarrowfill}}\,}}

\def\toPD{{\scriptstyle\,\buildrel{\P(\D)}%
\over{\hbox to 0.8cm{\rightarrowfill}}\,}}

\def\toDD{{\scriptstyle\,\buildrel{\D2}%
\over{\hbox to 0.6cm{\rightarrowfill}}\,}}

\def\car{\vrule height 4pt depth 0pt width 4pt}

\begin{document}

\date{}

\title{Natural decomposition of processes\\ and weak Dirichlet processes}

\maketitle

\centerline{\bf Fran\c cois Coquet$^*$, Adam
Jakubowski$^{**,}$\footnote{Supported in part by Komitet Bada\'n
Naukowych under Grant No 1 P03A 022 26 and completed while the
author was visiting Universit\'e de Rennes I.}}

\centerline{\bf Jean M\'emin$^{***}$ and Leszek S\l
omi\'nski$^{**,}$\footnote{Supported in part by Komitet Bada\'n
Naukowych under Grant No 1 P03A 022 26 and completed while the
author was visiting Universit\'e de Rennes I.}}

\medskip
\centerline{\it $^*$LMAH, Universit\'e du Havre}

\centerline{\it 25 rue Philippe Lebon,
76063 Le Havre Cedex, France}

\centerline{\it $^{**}$Faculty of Mathematics and Computer Science,}

 \centerline{\it N. Copernicus University, ul. Chopina, 87-100 Toru\'n,
Poland}

\centerline{\it $^{***}$IRMAR, Universit\'e de Rennes 1}

\centerline{\it Campus de Beaulieu,
35042 Rennes Cedex, France}




\bigskip

\begin{center}{\small{\bf Abstract}}
\end{center}
{ \small A class of stochastic processes, called ``weak Dirichlet
processes", is introduced and its properties are investigated in
detail. This class is much larger than the class of Dirichlet
processes. It is closed under $C^1$-transformations and under
absolutely continuous change of measure. If a weak Dirichlet
process has finite energy, as defined by Graversen and Rao, its
Doob-Meyer type decomposition is unique. The developed methods
have been applied to a study of generalized martingale
convolutions.}

\medskip

\medskip

{\it Mathematics Subject Classification:} 60G48, 60H05

\medskip

{\it Keywords:} weak Dirichlet processes, Dirichlet processes,
semimartingales, finite energy processes, quadratic variation,
mutual quadratic covariation, Ito's formula, generalized
martingale convolution.

\section{Introduction}

The quadratic variation of a stochastic process as well as the mutual
covariation
of two stochastic processes have been well-known for a long time to be at
the core
of the theory of stochastic integration, was it only because the quadratic
variation
appears
explicitely in Ito's formula for semimartingales. And indeed, every attempt to
 generalize Ito's calculus to a wider class of integrators (for instance
Dirichlet
 processes) or to functions less regular than ${\cal C}^2$ functions has to
deal
 with quadratic variations or covariations of the processes that appear. In
this
 aspect, the most enlightening work is perhaps F\"ollmer's paper (\cite{F850}).

On the other hand, it was proven by Graversen and Rao (\cite{GR}) that the
existence of a Doob-Meyer type decompositon for a process $X$ is narrowly
 linked to the fact that $X$ has a finite energy, which is a somewhat weaker
 assumption than the existence of a quadratic variation for $X$. The most
 well-known class of processes with finite energy (beyond the class of
semimartingales) is the class of Dirichlet processes. A larger class has
 been recently introduced by Errami and Russo (\cite{ER1}) under the name
 ``weak Dirichlet processes". The present paper explores some desirable
 properties of such processes. Although our definition of a quadratic
 variation is different from Errami and Russo's one -it is in a way more
classical- at any rate both
coincide as far as semimartingales are concerned. The other noticeable
difference is that throughout the
paper we deal with non continuous processes.

In part 2., we give an as explicit as possible link between
quadratic variation, energy, Dirichlet processes, weak Dirichlet
processes, and ``natural" (that is, ``Doob-Meyer type")
decomposition.

Part 3. is devoted to prove that any ${\cal C}^1$ function of a weak
Dirichlet process is again a weak
Dirichlet process. We are able to give an explicit Ito-type formula for
${\cal C}^2$ transformations,
but we could only find an explicit formula for the martingale part in the
general case.

Part 4. which is closest to Errami and Russo's work mentionned above, deals
with processes
 $X$ that may be written $\displaystyle X_t=\int_0^tG(t,s)dL_s$ where $L$
is a quasileft continuous
-but not necessarily continuous- square-integrable martingale and $G$ is a
deterministic function. We
give two sets of hypotheses under which $X$ is a weak Dirichlet process,
and also give its natural
decomposition. This section is illustrated through 3 examples, last one
giving additionnaly a formula
of Fubini type.

At last, we joined as an appendix some counter-examples related to
quadratic variation or to
 regularity of paths of processes. Although such examples may be
well-known, we could not find any in
the litterature, and we hope that they can enlight some of the technical
problems we are confronted
with here and there in the paper.

\section{Basic notations and results about processes with finite energy and
weak Dirichlet processes}

\medskip

In what follows, we are given a probability space $(\O, {\cal G}, \P )$.

We also fix a positive real number $T$. Unless otherwise stated, every process
or filtration will be indexed by $t\in [0, T]$. A filtration
$(\F_t)_{t\leq T}$
is denoted by $\F $. All filtrations  are assumed to be right-continuous
and defined on  $(\O, {\cal G}, \P )$ with
$\F_T \subset {\cal G}$.

We are also given a refining sequence ${D_n}$ of subdivisions of $[0,T]$ whose
mesh goes to 0 when $n\to\infty$. For every $n$, $D_n=\{0=t_0^n,t_1^n,\dots,
t_{N(n)}^n=T\}$.

We work with processes with a.s. right-continuous trajectories
 with left limits (such a process is called c\`adl\`ag), null in 0 and,
unless otherwise stated, admitting
a finite energy in the  sense defined below following  Graversen and Rao
(\cite{GR}):

\begin{definition} We say that $X$ is a process
{\it of finite energy} if

\begin{equation}
\label{EF}
\sup_nE\left[\sum_{t_i^n\in D_n}(X_{t_i^n}-X_{t_{i-1}^n})^2\right]<+\infty.
\end{equation}

This ``sup" will be denoted ${\cal E}n(X)$.

\end{definition}

\medskip

Of course, if $X$ has a finite energy,
${\displaystyle |X_t|^2}$ is integrable for every $t\leq T$ and
also $\sum_{s\leq T}\Delta X_s^2$ is integrable.

\medskip

We recall Graversen and Rao's main result in \cite{GR}:
\begin{theorem}\label{GRth}
If $X$ is a process with finite energy, then we can write $X$ as a sum
$X=M+A$,
where M is a square-integrable martingale and $A$ is a predictable process such
that there exists a subsequence $(D_{n_j})$ of $(D_n)$ satisfying
\begin{equation}
\label{M0}
E\left[\sum_{t_i^{n_j}\in D_{n_j},t_i^{n_j}\leq
t}(A_{t_i^{n_j}}-A_{t_{i-1}^{n_j
}})(N_{t_i^{n_j}}-N_{t_{i-1}^{n_j}})\right] \longrightarrow 0
\end{equation}
as $j\to \infty$ for all square integrable martingale $N$.

If, moreover, $X=M'+A'$ is any other such decomposition, the process $A-A'$ is
a continuous martingale.

At last, if we write
\begin{equation}\label{Mn}
M^n_t=  \sum_{t^n_i\in D_n, t^n_{i}\leq
t}\left[X_{t^n_{i}}-E[X_{t^n_{i}}/\F_{t^n_{i-1}}]
\right],
\end{equation}
then for all $t\in\bigcup_n D_n$, $M_t$ is the weak limit in $\sigma({\bf L}^2,
{\bf L}^2)$
of the sequence, $(M^{n_j}_t)$.
\end{theorem}

\bigskip

Such a decomposition of $X$ is a Doob-Meyer type decomposition:
the predictable process $A$ with the property of convergence (2)
is a ``natural" process.  In  this section we will discuss the
case when the decomposition $X=M+A$ in Theorem \ref{GRth} is
unique. Such Doob-Meyer decomposition will be called the natural
decomposition of $X$.

\medskip

We will use
the  notion of weak Dirichlet process introduced by Errami and Russo
(\cite{ER1}) in a slighty
different context.

\begin{definition}\label{wdp} We say that $X$ is a weak Dirichlet process
if it admits a
decomposition $X=M+A$, where $M$ is a
local martingale and $A$ is a predictable process such as $[A,N]=0$ for all
continuous local martingale $N$.
\end{definition}

In the above definition and in the sequel  we use the notion of  mutual
covariation and quadratic variation  in the following sense taken from
\cite{CMS}.

\begin{definition}{\bf (i)}
 Processes $X$ and $Y$ admit a quadratic (mutual) covariation along $(D_n)$ if
there exists a c\`adl\`ag process denoted $[X,Y]$ with for every $t\leq T$
$$[X,Y]_t=[X,Y]_t^c+\sum_{s\leq t}\Delta X_s\Delta Y_s$$
 and
\begin{equation}
\label{CVQ}
S^n(X,Y)_t:=\sum_{t_i^n\in D_n,t_i^n\leq
t}(X_{t_{i+1}^n}-X_{t_{i}^n})(Y_{t_{i+1}^n}-Y_{t_{i}^n}
)\toP [X,Y]_t\,\,\,\,as\,\, n\to\infty.
\end{equation}

{\bf (ii)} The process $X$ admits a quadratic variation along $(D_n)$ if
there exists a
c\`adl\`ag process denoted $[X,X]$ with for every $t\leq T$
$$[X,X]_t=[X,X]_t^c+
\sum_{s\leq t}\Delta X^2_s$$
 and
\begin{equation}
\label{VQ}
S^n(X,X)_t:=\sum_{t_i^n\in D_n,t_i^n\leq t}(X_{t_{i+1}^n}-X_{t_{i}^n})^2\toP
[X,X]_t\,\,\,\,as\,\, n\to\infty.
\end{equation}

\end{definition}

\medskip

\begin{remark} {\rm {\bf (i)} The decomposition
$X=M+A$ of a weak Dirichlet process is unique.

To see this suppose  that we have decompositions $X=M+A=M'+A'$ with  $A$
and  $A'$   predictable and verifying  $[A,N]=[A,N']=0$ for every
continuous local martingale  $N$. Then  $A-A'$  is a predictable local
martingale, hence a continuous local martingale. Then

$[A-A']_T=[A-A',A]_T-[A-A',A']_T=0$

and we deduce that $A=A'$.

{\bf (ii)} A weak Dirichlet process $X$ need not admit a quadratic
variation. We know only that for every continuous martingale  $N$ there
exists the covariation $[X,N]$.

{\bf (iii)} Of course, in general, a decomposition $X=M+A$ with a martingale
$M$ and a predictable process $A$, does not imply that $[A,N]=0$ for every
continuous martingale $N$, when $[A,N]$ exists. For example, take $A$ as a
continuous martingale and $N=A$.
}\end{remark}

The class of weak Dirichlet processses is much larger than the class of
Dirichlet processes. We recall:

\begin{definition}
A Dirichlet process is the sum of a local martingale and a continuous
process whose quadratic variation is identically zero.
\end{definition}

\begin{remark}{\rm
 Note that a Dirichlet process admits  a quadratic variation, which
is equal to the quadratic variation of its martingale part.
Our definition of a quadratic variation, which follows F\"ollmer's one
in \cite{F850} is weaker than the definition in \cite{F851}, and slightly
different from Russo and Vallois' one in \cite{RV}. However the three of them
coincide as far as semimartingales are concerned, and a Dirichlet process
according to the definition in \cite{F851} is also Dirichlet according to the
two other ones.}
\end{remark}
\medskip

The following notion of pre-quadratic variation is weaker than the quadratic
variation one; however the two notions coincide under stronger assumptions,
as will be seen below.

\begin{definition}
A process $X$ (not necessarily c\`adl\`ag) admits a pre-quadratic variation
along $(D_n)$ if there exists an increasing process denoted $S(X,X)$ with for
every $t\leq T$
\begin{equation}
\label{PVQ}
S^n(X,X)_t:=\sum_{t_i^n\in D_n,t_i^n\leq t}(X_{t_{i+1}^n}-X_{t_{i}^n})^2\toP
S(X,X)_t\,\,\,\,as\,\, n\to\infty.
\end{equation}
\end{definition}

\begin{remark}
{\rm We can find examples of continuous processes $X$ such that $S(X,X)$
is defined but not continuous (see example in Annex), hence $X$ does not
admit a quadratic variation.
}
\end{remark}

For every $t\leq T$ denoting $\pi_{t}$ any subdivision of $[0,t]$,
we consider the sum

$$S^{\pi_t}(X,X):=\sum_{{t_i\in \pi_t}, i>0}(X_{t_{i}}-X_{t_{i-1}})^2$$

\begin{proposition}
If a c\`adl\`ag process $X$ admits a pre-quadratic variation $S$ with the
following property $(S)$:

$S(X,X)$ is right continuous and for every $t\leq T$,
$S^{\pi_t}(X,X)\toP S(X,X)_t$
 as the mesh of $\pi_t$ tends to $0$.

Then, $S(X,X)$ is the quadratic variation of $X$ along any sequence $(D_n)$ of
subdivisions of $[0,T]$, whose mesh tends to $0$.
\end{proposition}

This result is proved in (\cite{J}), Lemme (3.11).

\begin{remark}
{\rm {\bf (i)} The class of Dirichlet processes is larger than
the space ${\bf H^2}$ of
semimartingales. Every continuous function admitting a quadratic variation
equal to zero is a deterministic Dirichlet process.

{\bf(ii)} Every continuous function is a deterministic weak
Dirichlet process:

Actually, let us consider a bounded continuous martingale $N$ nul in $0$,
we have
\begin{eqnarray*}
E(|\sum_{t_i^n\in D_n}(f(t_{i+1}^n)-f(t_i^n))
(N_{t_{i+1}^n}-N_{t_{i}^n})|^2)
&\leq \sup_{t_i^n}(f(t_{i+1}^n)-f(t_i^n))^2E(N_T)^2.
\end{eqnarray*}
and, from continuity of $f$  this last term tends to $0$ when $n\rightarrow
\infty $.

\medskip

 We give in Section 4
nondeterministic examples of weak Dirichlet processes, which are not ordinary
Dirichlet processes.
}

\end{remark}

\begin{remark}
{\rm The family of processes with finite energy is clearly stable under
addition, however we do not know if this stability holds for the family of
processes
admitting a quadratic variation. Of course this is true for the family of
Dirichlet processes. }
\end{remark}

\begin{theorem} Assume $X$ is a process with finite energy. The following
three conditions are equivalent:

\begin{description}

\item[{\bf(i)}]
$X$ is a weak Dirichlet process,

\item[{\bf(ii)}] for every
continuous local martingale $N$,
the quadratic covariation $[X,N]$ is well-defined,

\item[{\bf(iii)}] for every locally square integrable
martingale $N$, the quadratic covariation $[X,N]$ is well-defined.

\end{description}

In this case the decomposition
$X=M+A$ is unique and it is the natural decomposition expressed in Theorem
\ref{GRth}.
\end{theorem}

{\it Proof}: $(i)\Rightarrow(iii)$ Let us write $X=M+A$ as in Definition
\ref{wdp} and consider the
decomposition
$N=N^c+N^d$, where  $N^c$  is the continuous and $N^d$  purely
discontinuous part of $N$. By the definition of a weak Dirichlet process
the covariation $[X,N^c]$ is
well-defined. For proving the existence of $[X,N^d]$ we use
the following lemma.

\medskip

\begin{lemma}\label{scs} Assume that $X$ has a finite energy, and that $N$ is a
locally square integrable martingale which is the compensated sum of its
jumps. Then $X$ and $N$ admit a covariation such that
\begin{equation}
[X,N]_t=\sum_{s\leq t}\Delta X_s\Delta N_s.
\end{equation}
\end{lemma}

{\it Proof of Lemma \ref{scs}}: By using a localizing sequence of stopping
times,
one can assume that
$N$ is a square integrable martingale. One can find a sequence $(N^p)_p$ of
martingales with finite
variation and only a finite number of jumps, such that $N^p\to N$ in ${\bf
H^2}$.

We have then, for fixed $p$,
\begin{equation}
\sum_{t_i^n\in D_n,t_i^n\leq
t}(X_{t_{i+1}^n}-X_{t_{i}^n})(N^p_{t_{i+1}^n}-N^p_{t_{i}^n}
)\toP \sum_{s\leq t}\Delta X_s\Delta N^p_s.
\end{equation}
as $n\to\infty$.

\medskip

On the other hand,
\begin{eqnarray*}
&&E\left|\sum_{t_i^n\in D_n,t_i^n\leq t}(X_{t_{i+1}^n}-X_{t_{i}^n})
(N^p_{t_{i+1}^n}-N^p_{t_{i}^n})-
\sum_{t_i^n\in D_n,t_i^n\leq t}(X_{t_{i+1}^n}-X_{t_{i}^n})(N_{t_{i+1}^n}-
N_{t_{i}^n})\right|\\
&&\qquad\leq E|\Bigl(\sum_{t_i^n\in D_n,t_i^n\leq t}
(X_{t_{i+1}^n}-X_{t_{i}^n})^2\Bigr)^{1/2}\\
&&\qquad\qquad\qquad\times\Bigl(\sum_{t_i^n\in D_n,t_i^n\leq t}\bigl
((N_{t_{i+1}^n}-N^p_{t_{i+1}^n})
-(N_{t_{i}^n}-N^p_{t_{i}^n})\bigr)^2\Bigr)^{1/2}|\\
&&\qquad\leq ({\cal E}n(X))^{1/2}E\Bigl([N-N^p,N-N^p]_t\Bigr)^{1/2}
\end{eqnarray*}

which goes to 0 as $p\to\infty$ since $[N-N^p,N-N^p]_t\toP 0$.

\medskip

At last,
\begin{eqnarray*}
& &E\Bigl|\sum_{s\leq t}\Delta X_s(\Delta N_s-\Delta N^p_s)\Bigr|\leq
E\Bigl(\bigl(
\sum_{s\leq t}\Delta X_s^2\bigr)^{1/2}\bigl(\sum_{s\leq t}(\Delta N_s-\Delta
N^p_s)^2\bigr)^{1/2}\\
&& \qquad\qquad\qquad\qquad
\leq ({\cal E}n(X))^{1/2}E\left[[N-N^p,N-N^p]_t\right]^{1/2}\cr
\end{eqnarray*}
which goes to zero as $p$ goes to infinity, hence $\sum_{s\leq t}\Delta X_s
\Delta N^p_s$ converges in $L^1$ to $\sum_{s\leq t}\Delta X_s\Delta N_s$.

These three convergences give the lemma.
\car\car

\bigskip

$(iii)\Rightarrow(ii)$ is obvious.

\medskip

$(ii)\Rightarrow(i)$ Let $X=M+A$ be a decomposition from  Theorem \ref{GRth}
 and let $N$ be a continuous local martingale. Define
$T_p=\inf\{ t:\vert N_t\vert \geq p\} $ then $(T_p)$ is a
localizing sequence of stopping times. We will prove that for every $p$,
$[A,N^{T_p}]=0$, which implies that also $[A,N]=0$.

By hypothesis we have the convergence $S^n(X,N^{T_p})_t\toP [X,N^{T_p}]_t$.

Hence we deduce that also $S^n(A,N^{T_p})_t\toP [A,N^{T_p}]_t$. We will
show that in fact
\begin{equation}\label{eqnew}
S^n(A,N^{T_p})_{T}\rightarrow [A,N^{T_p}]_{T}\quad
in \,\,L^1.\end{equation}
To see this, it is sufficient to check uniform integrability of
the sequence $\{S^n(A,N^{T_p})_{T}\}$.
Writing $\vert S^n(X,N^{T_p})_{T}\vert\leq S^n(X,X)^{1/2}_{T}
S^n(N^{T_p},N^{T_p})^{1/2}
_{T}$, and using H\"older inequality, we get:

\[E\vert S^n(X,N^{T_p})_{T}\vert ^{4/3}
\leq E[S^n(X,X)_{T}]E[S^n(N^{T_p},N^{T_p})_{T}^2].\]

Similarly we prove that also

$E\vert S^n(M,N^{T_p})_{T}\vert^{4/3}<+\infty$.

As a consequence,
we deduce uniform integrability of $\{(S^n(A,N^{T_p})_{T}\}$ and
(\ref{eqnew}) holds true. Therefore, in particular
\[E[S^n(A,N^{T_p})_{T}]\rightarrow E[A,N^{T_p}]_{T}\]
and due to  (\ref{M0} )
\begin{equation}\label{eqnew1}
E[A,N^{T_p}]_{T}=0.\end{equation}

Note that the process $[A,N^{T_p}]$
 has a finite variation ; moreover, since $N$
is continuous, $[A,N^{T_p}]$ is also a continuous process. Therefore,
to get $[A,N^{T_p}]=0$
it is  sufficient to prove that $[A,N^{T_p}]$ is a local martingale

Let us consider a bounded stopping time $\tau\leq T $,
the same arguments as  above give the convergence:

$E[S^n(A,N^{T_p\wedge \tau })_{T}]
\rightarrow E[A,N^{T_p\wedge \tau }]_{T}$

and (\ref{M0} ) gives

\[E[A,N^{T_p }]_{\tau}=E[A,N^{T_p\wedge \tau }]_{T }=0.\]

Hence, it follows easily that the stopped process $[A,N]^{T_p}$ is a martingale
and $[A,N^{T_p }]=0$. Since $P(T_p=T)\uparrow 1$, the proof of the last
implication is completed.

Finally, note that the uniqueness of the decomposition in Theorem \ref{GRth}
is an easy consequence of  the fact that  $X$ is a weak Dirichlet process.
\hfill{\car} \\[2mm]
\car

We get immediately the following
\begin{corollary} Let us consider a weak Dirichlet process $X$

{\bf (i)} If $Q$ is a probability measure absolutely continuous with respect to
$\P $, then $X$ is a $Q$ weak Dirichlet process.

{\bf (ii)} For an $a>0$ we define $\hat X= \sum_{s\leq .}\Delta
X_s1_{\Delta |X_s|
>a} $, then $X-\hat X$ is a weak Dirichlet process.
\end{corollary}

Now,  we will consider processes with finite energy $X$
admitting additionally a quadratic variation $[X,X]$. Then of course
$E[X,X]_T<\infty$.

\begin{theorem}\label{thm2} Assume $X$ is a weak Dirichlet
process with finite energy, admitting a
quadratic  variation process.
\begin{description}
\item[{\bf(i)}] In the natural decomposition $X=M+A$, $M$ is
a square integrable martingale and $A$ has an integrable quadratic variation.

\item[{\bf(ii)}]
The natural decomposition is minimal in the following sense.
If $X=M'+A'$ is another decomposition with a local martingale
$M'$ and a predictable process $A'$

then  $[A'A']$ is well defined and:
$$[A',A']=[M-M',M-M']+[A,A]. $$
\end{description}

\end{theorem}

{\it Proof}: $(i)$ To begin with, we notice that

$$E[\sum_{s\leq T}\Delta A_s^2]<\infty; $$

actually,
for every predictable stopping time $S$, $\Delta A_S=E[\Delta X_S\vert
\F_S-]$, hence
$\displaystyle E[\sum_{s\leq T}\Delta A_s^2]\leq E[\sum_{s\leq T}\Delta
X_s^2]<\infty.$
It follows that
$\displaystyle E[\sum_{s\leq T}\Delta M_s^2]<\infty $ and $M$ is a locally
square integrable
martingale.

Let us consider the decomposition $M=M^c+M^d$ where $M^c$ is the continuous
part of $M$
and $M^d$ its purely discontinuous part. Writing
$[M^d,A]=[M^d,X]-[M^d,M^d]$, we
get the existence of $[M^d,A]$. Using the property of decomposition
$X=M+A$, we have
$[M,A]=[M^d,A]$. From Lemma \ref{scs} we deduce:

$$[M,A]=\sum_{s\leq \cdot}\Delta M_s\Delta X_s-\sum_{s\leq \cdot}
\Delta M_s^2=\sum_{s\leq \cdot}\Delta M_s
\Delta A_s.$$

Now, by the definition of quadratic variation of $X$ and $M$
one gets the existence of $[A,A]$:

$$[A,A]=[X,X]-2[M,A]-[M,M]. $$

Finally,
$$E([M,M]_T+[A,A]_T)\leq E[X,X]_T+2E[\sum_{s\leq T}\Delta M_s^2]^{1/2}
E[\sum_{s\leq T}\Delta A_s^2]^{1/2}<\infty .$$

(ii) Since $A'=M+A-M'$, by linearity $[A',A']$ is well defined and we
can write:

$$[A',A']=[M+A-M',M+A-M']=[M-M',M-M']+[A,A]-2[A,M-M']. $$

But, $M-M'$ is a continuous local martingale and as $A$ is taken from the
natural decomposition of $X$, we
get: $[A,M-M']=0$, hence the desired result.

\hfill{\car}

\bigskip

\section{Stability of weak Dirichlet processes under $C^1$ transformations}

Assume $X$ is a process of finite energy. Let us denote by
$\mu $ the jump measure of $X$. Then $\displaystyle\sum_{s\leq .}\Delta
X_s^2$
can be written $\displaystyle\int_0^.\int_{\R -\{ 0\} }x^2\mu (ds,dx)$. Since
$\displaystyle E\sum_{s\leq T}\Delta X^2_s<\infty $, $X$ admits a L\'evy system
$\nu $ which is the predictable compensator of $\mu $; then the
predictable
increasing process $\displaystyle\int_0^.\int_{\R -\{ 0\} }x^2\nu (ds,dx)$
is well
defined and
$$E\sum_{s\leq T}\Delta X^2_s=E[\int_0^T\int_{\R -\{ 0\} }x^2\mu (ds,dx)]=
E[\int_0^T\int_{\R -\{ 0\} }x^2\nu (ds,dx)]<\infty.$$

\medskip

We begin with $C^2$ stability:

\medskip

\begin{theorem}\label{C2}  Let $X=M+A$ be a weak Dirichlet process of
finite energy and $F$ a $C^2$-real valued
function with  bounded derivatives $f$ and $f'$. Then the process
$(F(X_t)_{t\geq 0})$ is a
weak Dirichlet process of finite energy
and the decomposition $F(X)= Y+\Gamma $ holds with the martingale part

\begin{eqnarray*}& &Y_t= F(0)+\int_0^tf(X_{s-})dM_s\\
&&\qquad+\int_0^t\int_{\R }\Big(F(X_{s-}+x)-F(X_{s-})-xf(X_{s-})\Big )
(\mu -\nu )(ds, dx)
\end{eqnarray*}
and the predictable part
\begin{eqnarray*} \lefteqn{\Gamma
_t=\int_0^tf(X_{s-})dA_s-1/2\sum_{s\leq t}f'(X_{s-})(\Delta A_s)^2
+1/2\int_0^tf'(X_s)d[M,M]^c_s}\cr
&\qquad+\int_0^t\int_{\R
}\Big(F(X_{s-}+x)-F(X_{s-})-xf(X_{s-})\Big)\nu (ds ,dx),
\end{eqnarray*}
where  $\displaystyle(S)\int_0^. f(X_{s-})dA_s$ is well defined as a
limit in probability of Riemann sums. More precisely for every $t$
$$\sum_{t^n_i\in D_n, t^n_i\leq
t}(f(X_{t^n_i})(A_{t^n_{i+1}}-A_{t^n_i})+1/2f'(X_{t^n_i})
(A_{t^n_{i+1}}-A_{t^n_i})^2)
\toP (S)\int_0^t
f(X_{s-})dA_s.$$
\end{theorem}

{\it Proof:} Fix  $t>0$.
We use arguments from the paper \cite{F850} by F\"ollmer. For $\epsilon>0$ we
define  $J(1)=\{s\leq
t;|\Delta X_s|>\epsilon\}$. In the following, the elements of $D_n$ are,
for short, written $t_{i}$ instead of $t^n_i$. Then

$ \sum_{t_i\in D_n, t_i\leq
t}(F(X_{t_{i+1}})-F(X_{t_i}))=\sum_1(F(X_{t_{i+1}})-F(X_{t_i}))
+\sum_2(F(X_{t_{i+1}})-F(X_{t_i})),$

where $\displaystyle\sum_1$ denotes the sum (depending on $\omega\in\Omega$)
of  $t_i\in D_n, t_i\leq
t$ such that  $(t_i,t_{i+1}]\cap J(1)\neq\emptyset$ and
$\displaystyle\sum_2$ the sum on the other $t_i$'s.

Then by Taylor's formula
\begin{eqnarray*}
& &\sum_2(F(X_{t_{i+1}})-F(X_{t_i}))
=\sum_2 f(X_{t_i})(X_{t_{i+1}}-X_{t_i})\\
& &\nonumber\qquad\qquad+1/2\sum_2 f'(X_{t_i})(X_{t_{i+1}}-X_{t_i})^2+\sum_2
r_2(X_{t_i},X_{t_{i+1}}),
\end{eqnarray*}
where $r_2(X_{t_i},X_{t_{i+1}})=C^{\epsilon}_i(X_{t_{i+1}}-X_{t_i})^2$ with
$\max_2
|C^{\epsilon}_i|\leq ||f'||$ and
\begin{equation}\label{eqnew2}
\lim_{\epsilon\downarrow0}\limsup_{n\to \infty} P(\max_2
|C_i^{\epsilon}|>\delta)=0,\quad \delta>0.\end{equation}

Hence
\begin{eqnarray*}
& &\sum_{t_i\in D_n, t_i\leq
t}(F(X_{t_{i+1}})-F(X_{t_i}))
=\sum_{t_i\in D_n, t_i\leq
t} f(X_{t_i})(X_{t_{i+1}}-X_{t_i})\\
& &\nonumber+1/2\sum_{t_i\in D_n, t_i\leq
t} f'(X_{t_i})(X_{t_{i+1}}-X_{t_i})^2+\sum_2
r_2(X_{t_i},X_{t_{i+1}})\\
& &
-\sum_1\{F(X_{t_{i+1}})-F(X_{t_i})-f(X_{t_i})(X_{t_{i+1}}-X_{t_i})
-1/2f'(X_{t_i})(X_{t_{i+1}}-X_{t_i})^2\}\\
& &=\sum_{t_i\in D_n, t_i\leq
t} f(X_{t_i})(M_{t_{i+1}}-M_{t_i})+\sum_{t_i\in D_n, t_i\leq t}
f(X_{t_i})(A_{t_{i+1}}-A_{t_i})\\
& & +1/2\sum_{t_i\in D_n, t_i\leq
t} f'(X_{t_i})(A_{t_{i+1}}-A_{t_i})^2+1/2\sum_{t_i\in D_n, t_i\leq
t}
f'(X_{t_i})(M_{t_{i+1}}-M_{t_i})^2\\
& &+\sum_{t_i\in D_n, t_i\leq
t} f'(X_{t_i})(M_{t_{i+1}}-M_{t_i})(A_{t_{i+1}}-A_{t_i})
+ \sum_2r_2(X_{t_i},X_{t_{i+1}})\\
& &
-\sum_1\{F(X_{t_{i+1}})-F(X_{t_i})-f(X_{t_i})(X_{t_{i+1}}-X_{t_i})
-1/2f'(X_{t_i})(X_{t_{i+1}}-X_{t_i})^2\}\\
& & I^n_1+I^n_2+I^n_3+I^n_4+I^n_5+I^{n,\epsilon}_6
+I^{n,\epsilon}_7.\end{eqnarray*}

Now,  note that by the definition of a stochastic integral we have

$I^n_1\toP\int_0^tf(X_{s-})dM_s. $

\medskip

The following  simple lemma will be very useful in order to estimate the
other terms.

\begin{lemma}\label{lemnew}
Assume that c\`adl\`ag processes $X$ and $Y$ admit a quadratic
covariation $[X, Y]$ and that the sequence
$\{Var(S^n(X,Y))_T\}$ is bounded in probability.

To c\`adl\`ag processes $Z$ and $U$ we associate the sequences $\{Z^n\}$ and
$\{U^n\}$  of processes, where
$Z^n$ and $U^n$ are the respective discretizations of $Z$ and $U$ along
$D_n$; precisely $Z^n_t=
Z_{t^n_{i}}$, $U^n_t=
U_{t^n_{i}}$ when
$t\in [t^n_i, t^n_{i+1}[$. Then, for every continuous real function
$f,g$ and every $t$, holds the convergence:
$$\int_0^tf(Z^n_{s-})g(\Delta U^n_s)dS^n(X,Y)_s \toP
\int_0^tf(Z_{s-})g(\Delta U_s)d[X,Y]_s,$$
where these integrals are Stieltjes integrals with
respect to the
processes $S^n(X,Y)$ or $[X,Y]$.
\end{lemma}

{\it Proof of Lemma 3.1} From the proof of  \cite[Lemma 1.3]{CMS}
one can  deduce that

$\int_0^.f(Z^n_{s-})g(\Delta U^n_s)dS^n(X,Y)_s
\toP\int_0^.f(Z_{s-})g(\Delta U^n_s)d[X,Y]_s $

in the (so-called $J_1$)  Skorokhod topology  (see e.g.
\cite{JS}). Since for every $t$

$f(Z^n_{t-})g(\Delta U^n_t)\Delta S^n(X,Y)_t\toP f(Z_{t-})g(\Delta U_t)
\Delta [X,Y]_t, $

the desired
result  follows from properties
of  the Skorokhod topology $J_1$. \hfill{\car} \\[2mm]

It is clear by Lemma \ref{lemnew} that
we have the convergences;
\begin{eqnarray*}
& &I^n_4\toP1/2\int_0^tf'(X_{s-})d[M,M]_s,\\
& &I^n_5\toP\int_0^tf'(X_{s-})d[M^d,A]_s,
\end{eqnarray*}
where $M^d$ denotes purely discontinuous part of $M$.

Since  $X$ is a process with finite energy, by (\ref{eqnew2})

$
\lim_{\epsilon\downarrow0}\limsup_{n\to \infty} P(
|I^{n,\epsilon}_6|>\delta)=0,\quad \delta>0. $

We observe also that $P$-almost surely there exists the limit
\begin{eqnarray*}
& &\lim_{\epsilon\downarrow0}\lim_{n\to \infty}I^{n,\epsilon}_7=\sum_{s\leq
t}\{F(X_s)-F(X_{s-})-f(X_{s-})\Delta X_s-1/2f'(X_{s-})\Delta X_s^2\}\\
& &\qquad=\sum_{s\leq
t}\{F(X_s)-F(X_{s-})-f(X_{s-})\Delta X_s\}\\
& &-1/2\sum_{s\leq
t}\{f'(X_{s-})\Delta M_s^2+\sum_{s\leq
t}f'(X_{s-})\Delta A_s^2+2\sum_{s\leq
t}f'(X_{s-})\Delta M^d_s\Delta A_s\}.
\end{eqnarray*}
On the other hand it is obvious that $P$-almost surely

$\sum_{t_i\in D_n, t_i\leq
t}(F(X_{t_{i+1}})-F(X_{t_i}))\rightarrow F(X_t)-F(0) $

and putting together all convergences, we deduce that $\{I^n_2+I^n_3\}$ is
converging in probability and the limit we denote as
$(S)\int_0^tf(X_{s-})dA_s$. Observe also that

 $\int_0^tf'(X_{s-})d[M,M]^c_s=\int_0^tf'(X_{s-})d[M,M]_s-\sum_{s\leq
t}f'(X_{s-}) \Delta M_s^2 $

and
$$\int_0^tf'(X_{s-})d[M^d,A]_s=\sum_{s\leq t}f(X_{s-})\Delta M_s\Delta A_s.$$
As a consequence we obtain the formula
\begin{eqnarray*}
& & F(X_t)= F(0)+\int_0^tf(X_{s-})dM_s +(S)\int_0^tf(X_{s-})dA_s\\
& &\qquad\qquad+1/2\int_0^tf'(X_{s-})d[M,M]^c_s-1/2\sum_{s\leq
t}f'(X_{s-})\Delta A_s^2\\
& &\qquad\qquad+\sum_{s\leq t}\{ F(X_s)-F(X_{s-})-\Delta
X_sf(X_{s-})\}.
\end{eqnarray*}

Now, writing
\begin{eqnarray*}
&&\sum_{s\leq t}\{ F(X_s)-F(X_{s-})-\Delta
X_sf(X_{s-})\}\\
&&\qquad= \int_0^t\int_{\R -\{ 0\} }(F(X_{s-}+x)-F(X_{s-})-xf(X_{s-}))
\mu (ds,dx)
\end{eqnarray*}
and using the basic inequalities
$$\vert F(y+x)-F(y)-xf(y)\vert
\leq \Vert f'\Vert x^2$$
and
$$\vert F(y+x)-F(y)-xf(y)\vert \leq 2\Vert f\Vert \vert x\vert, $$
we get the decomposition
\begin{eqnarray*}
&&\sum_{s\leq .}\{F(X_s)-F(X_{s-})-\Delta
X_sf(X_{s-})\}\\
&&\qquad= \int_0^.\int_{\R -\{ 0\} }(F(X_{s-}+x)-F(X_{s-})-xf(X_{s-}))
(\mu -\nu )(ds,dx)\\
&&\qquad\qquad+\int_0^.\int_{\R -\{ 0\} }(F(X_{s-}+x)-F(X_{s-})-xf(X_{s-}))
\nu (ds,dx)
\end{eqnarray*}
where
$$\int_0^.\int_{\R -\{ 0\} }(F(X_{s-}+x)-F(X_{s-})-xf(X_{s-})
(\mu -\nu )(ds,dx)$$
is a square integrable purely discontinuous
martingale, which we will denote by $L$ and
$$\int_0^.\int_{\R -\{
0\} }(F(X_{s-}+x)-F(X_{s-})-xf(X_{s-}))\nu (ds,dx)$$
is an increasing predictable square integrable process.

Then we get the decomposition $F(X_.)=F(0)+Y_.+\Gamma _.$, as written in the
statement of Theorem 3.1.

It remains to prove that, for every  continuous local martingale
$N$, holds the equality: $[\Gamma ,N]=0$.

First note that
\begin{eqnarray*}
& &\sum_{t_i\in D_n, t_i\leq t}
(\Gamma_{t_{i+1}}-\Gamma_{t_i})(N_{t_{i+1}}-N_{t_i})\\
& &\qquad\qquad=
\sum_{t_i\in D_n, t_i\leq t}
\int_{t_i}^{t_{i+1}}(f(X_{s-})-f(X_{t_i}))dM_s(N_{t_{i+1}}-N_{t_i})\\
& &\nonumber\qquad\qquad\quad
-\sum_{t_i\in D_n, t_i\leq t}
(L_{t_{i+1}}-L_{t_i})(N_{t_{i+1}}-N_{t_i})\\
& &\qquad\qquad\quad+\sum_{t_i\in D_n, t_i\leq t}
f(X_{t_i})(A_{t_{i+1}}-A_{t_i})(N_{t_{i+1}}-N_{t_i})\\
& &\qquad\qquad\quad\nonumber +1/2\sum_{t_i\in D_n, t_i\leq t}
 f'(X_{t_i})(A_{t_{i+1}}-A_{t_i})^2(N_{t_{i+1}}-N_{t_i})\\
& &\qquad\qquad\quad+1/2\sum_{t_i\in D_n, t_i\leq t}
f'(X_{t_i})(M_{t_{i+1}}-M_{t_i})^2(N_{t_{i+1}}-N_{t_i})\\
& &\qquad\qquad\quad\nonumber+
\sum_{t_i\in D_n, t_i\leq t}
 f'(X_{t_i})(M_{t_{i+1}}-M_{t_i})(A_{t_{i+1}}-A_{t_i})(N_{t_{i+1}}-N_{t_i})
\\&&\qquad\qquad\quad
+ \sum_2r_2(X_{t_i},X_{t_{i+1}})(N_{t_{i+1}}-N_{t_i})\\
& &\qquad\qquad\quad
-\sum_1\{F(X_{t_{i+1}})-F(X_{t_i})-f(X_{t_i})(X_{t_{i+1}}-X_{t_i})\\
& &\qquad\qquad\qquad\quad
-1/2f'(X_{t_i})(X_{t_{i+1}}-X_{t_i})^2\}(N_{t_{i+1}}-N_{t_i})\\
& &\qquad\qquad=\nonumber I^n_1+I^n_2+I^n_3+I^n_4+I^n_5+I^n_6+I^{n,\epsilon}_7
+I^{n,\epsilon}_8\end{eqnarray*}
Clearly,

\[|I^n_1|\leq(\sum_{t_i\in D_n, t_i\leq t}
(\int_{t_i}^{t_{i+1}}(f(X_{s-})-f(X_{t_i}))dM_s)^2)^{1/2}
(\sum_{t_i\in D_n, t_i\leq t}(N_{t_{i+1}}-N_{t_i})^2)^{1/2},\]

where by the definition of the stochastic integral
\begin{eqnarray*}
& &E\sum_{t_i\in D_n, t_i\leq t}
(\int_{t_i}^{t_{i+1}}(f(X_{s-})-f(X_{t_i}))dM_s)^2
\\
& &\qquad\qquad=E\sum_{t_i\in D_n, t_i\leq t}
\int_{t_i}^{t_{i+1}}(f(X_{s-})-f(X_{t_i}))^2d[M,M]_s\rightarrow0.\end{eqnarray*}
On the other hand by Lemma \ref{lemnew}
\begin{eqnarray*}
& &I^n_2\toP[L,N]_t=0,\\
& &I^n_3\toP\int_0^tf(X_{s-})d[A,N]_s=0,\\
& &I^n_4\toP1/2\int_0^tf'(X_{s-})\Delta A_sd[A,N]_s=0,\\
& &I^n_5\toP1/2\int_0^tf'(X_{s-})\Delta N_sd[M,M]_s=0,
\end{eqnarray*}
and

\[I^n_6\toP\int_0^tf'(X_{s-})\Delta M_sd[A,N]_s=0.\]

Finally, for every  $\epsilon>0$

\[I^{n,\epsilon}_6\toP0\qquad\mbox{\rm and}\qquad
I^{n,\epsilon}_7\rightarrow0,\,\,P-a.s. \]

and the proof of Theorem \ref{C2} is completed. \hfill{\car}

\begin{corollary}\label{corC2}  Let $X=M+A$ be a weak Dirichlet process of
finite energy admitting a quadratic variation  and $F$ be a
$C^2$-real valued
function with  bounded derivatives $f$ and $f'$. Then the process
$(F(X_t)_{t\geq 0})$ is a
weak Dirichlet process of finite energy admitting a qudratic variation
and the decomposition $F(X)= Y+\Gamma $ holds with the
martingale part
\begin{eqnarray*}& &Y_t= F(0)+\int_0^tf(X_{s-})dM_s\\
&&\qquad+\int_0^t\int_{\R }\Big(F(X_{s-}+x)-F(X_{s-})-xf(X_{s-})\Big )
(\mu -\nu )(ds, dx)
\end{eqnarray*}
and the predictable part
\begin{eqnarray*} \lefteqn{\Gamma
_t=\int_0^tf(X_{s-})dA_s
 +1/2\int_0^tf'(X_s)
d[X,X]^c_s}\cr
&\qquad+\int_0^t\int_{\R
}\Big(F(X_{s-}+x)-F(X_{s-})-xf(X_{s-})\Big)\nu
(ds ,dx).
\end{eqnarray*}
\end{corollary}

{\it Proof}: By Theorem \ref{thm2}(i), $A$ admits integrable quadratic
variation $[A,A]$ and due to Lemma \ref{lemnew}

\[
1/2\sum_{t_i\in D_n, t_i\leq
t}f(X_{t_i})(A_{t_{i+1}}-A_{t_i})^2\toP1/2\int_0^tf(X_{s-})d[A,A]_s. \]

Since
$[X,X]^c=[M,M]^c+[A,A]^c$,
it is clear that
\begin{eqnarray*}
& &1/2\int_0^tf(X_{s-})d[A,A]_s-1/2\sum_{s\leq
t}f(X_{s-})\Delta A^2_s+1/2\int_0^tf(X_{s-})d[M,M]_s\\
& &\qquad\qquad=1/2\int_0^tf(X_{s-})d([A,A]^c_s
+[M,M]^c)=1/2\int_0^tf(X_{s-})d[X,X]^c_s
\end{eqnarray*}
and the decomposition $F(X_t)=F(0)+Y_t+\Gamma_t$ in the statement of
Corollary 3.1 is a consequence of
Theorem \ref{C2}.

Finally, by the Theorem from \cite[page 144]{F850} we obtain that also $F(X)$
admits a quadratic variation, which completes the proof.
\hfill{\car}

\begin{theorem}\label{C1}  Let $X=M+A$ be a weak Dirichlet process of
finite energy
and $F$ a $C^1$-real valued
function with  bounded derivative $f$.

Then the process $(F(X_t)_{t\geq 0})$ is a
weak Dirichlet process of finite energy
and the decomposition $F(X)= Y+\Gamma $ holds with the
martingale part
\begin{eqnarray*}& &Y_t= F(0)+ \int_0^tf(X_{s-})dM_s\\
&&\qquad+\int_0^t\int_{\R }\Big(F(X_{s-}+x)-F(X_{s-})-xf(X_{s-})\Big )
(\mu -\nu )(ds, dx)
\end{eqnarray*}
\end{theorem}

\begin{remark}
This theorem has formally almost the same statement as Theorem \ref{C2}.
However, we have not here any explicit formula for $\Gamma$. The
delicate point here is the behaviour of the sum

$$\sum_{s\leq t}\big (F(X_s)-F(X_{s-})-\Delta
X_sf(X_{s-})\big )$$

which is not necessarily absolutely convergent and which
does not define a process with finite variation.
\end{remark}

{\it Proof of Theorem 3.2}: We consider a sequence $(F^p)_{p\in \N }$ of $C^2$
real functions
such that $\Vert F-F^p\Vert +\Vert f-f^p\Vert \rightarrow 0$, when
$p\rightarrow
\infty $. Using Theorem \ref{C2} we can write
$$F^p(X_t)=F^p(0)+Y^p_t+\Gamma ^p_t$$
where
$$Y^p_t=\int_0^tf^p(X_{s-})dM_s+L^p_t$$
with
$$L^p_t=\int_0^t\int_{\R -\{ 0\}
}\Big(F^p(X_{s-}+x)-F^p(X_{s-})-xf^p(X_{s-})\Big )
(\mu -\nu )(ds, dx).$$
The sequence $\{\int_0^.f^p(X_{s-})dM_s+L^p_.\}_{p\in \N }$ is a Cauchy
sequence
in the space $H^2$ of square integrable martingales, hence the limiting
martingale
exists and has the form $\int_0^.f(X_{s-})dM_s+L_.$:

Actually, for $p,q$ integers
\begin{eqnarray*}
& &\Vert Y^p-Y^q\Vert _{H^2}\leq
E(\int_0^T(f^p(X_{s-})-f^q(X_{s-}))^2d[M,M]_s)\\
& &\qquad\qquad\qquad
+E[\int_0^T\int_{\R -\{ 0\} }\Big(F^p(X_{s-}+x)-F^p(X_{s-})-xf^p(X_{s-})\\
&&\qquad\quad\qquad\qquad\qquad
-F^q(X_{s-}+x)+F^q(X_{s-})+xf^q(X_{s-})\Big )^2\nu (ds, dx)]\\
& &\qquad\qquad\leq \Vert F^p-F^q\Vert ^2E[[M,M]]+2\Vert f^p-f^q\Vert
^2E[[M,M]].
\end{eqnarray*}
Now we write:

\[\Gamma^p_t= F^p(X_t)-F^p(0)-\int_0^tf^p(X_{s-})dM_s-L^p_t.\]

Clearly the sequence of predictable processes $(\Gamma^p)$ converges uniformly
in probability and its limit (i.e. the process $\Gamma$)
has to be also predictable.

It remains to prove that $[\Gamma,N]=0$ for every continuous local
martingale $N$.

 Fix   $t$.
In the sequel use the notations from the proof of Theorem \ref{C2}.

By Taylor's formula

\[ \sum_2(F(X_{t_{i+1}})-F(X_{t_i}))
=\sum_2 f(X_{t_i})(X_{t_{i+1}}-X_{t_i})\\
+\sum_2
r_1(X_{t_i},X_{t_{i+1}}), \]

where $r_1(X_{t_i},X_{t_{i+1}})=C_i^{\epsilon}
(X_{t_{i+1}}-X_{t_i})$ satisfy
$|C_i^{\epsilon}|\leq ||f||$ and

\[\lim_{\epsilon\downarrow0}\limsup_{n\to \infty} P(\max_2
|C_i^{\epsilon}|>\delta)=0,\quad \delta>0.\]

Therefore,
\begin{eqnarray*}
& &\sum_{t_i\in D_n, t_i\leq t}
(\Gamma_{t_{i+1}}-\Gamma_{t_i})(N_{t_{i+1}}-N_{t_i})
\\& &\qquad=\sum_{t_i\in D_n, t_i\leq
t} \int_{t_i}^{t_{i+1}}f(X_{t_i})-f(X_{s-})dM_s(N_{t_{i+1}}-N_{t_i})
\\& &\qquad\quad-\sum_{t_i\in D_n, t_i\leq t}
(L_{t_{i+1}}-L_{t_i})(N_{t_{i+1}}-N_{t_i})\\
& &\qquad\quad+\sum_{t_i\in D_n, t_i\leq t}
f(X_{t_i})(A_{t_{i+1}}-A_{t_i})(N_{t_{i+1}}-N_{t_i})\\
& &\qquad\quad\nonumber
+ \sum_2r_1(X_{t_i},X_{t_{i+1}})(N_{t_{i+1}}-N_{t_i})\\
& &\quad\quad\nonumber
-\sum_1\{F(X_{t_{i+1}})-F(X_{t_i})-f(X_{t_i})
(X_{t_{i+1}}-X_{t_i})(N_{t_{i+1}}-N_{t_i})\}\\
& &\qquad=\nonumber I^n_1+I^n_2+I^n_3+I^{n,\epsilon}_4
+I^{n,\epsilon}_5.\end{eqnarray*}
Clearly, first three sums tend to $0$ analogously to the proof of Theorem
\ref{C2}. Next,

\[\lim_{\epsilon\downarrow0}\limsup_{n\to\infty}P(|I^{n,\epsilon}_4|>\delta)=0,\
quad
\delta>0\]

and for every $\epsilon>0$

\[ I^{n,\epsilon}_5\rightarrow0,\,\,P-a.s. \]

which, completes  the proof of Theorem \ref{C1}. \hfill{\car}

\begin{corollary}  Let $X=M+A$ be a weak Dirichlet process of finite energy
admitting a quadratic variation process and $F$ a
$C^1$-real valued
function with  bounded derivative $f$. Then the process
$(F(X_t)_{t\geq 0})$ is a
weak Dirichlet process of finite energy admitting a quadratic variation
and the decomposition $F(X)= Y+\Gamma $ holds with the
martingale part
\begin{eqnarray*}& &Y_t= F(0)+ \int_0^tf(X_{s-})dM_s\\
&&\qquad+\int_0^t\int_{\R }\Big(F(X_{s-}+x)-F(X_{s-})-xf(X_{s-})\Big )
(\mu -\nu )(ds, dx)
\end{eqnarray*}
The quadratic variation process of $F(X_t)_{t}$
is given by

\[[F(X), F(X)]_t=\int_0^t (f(X_s))^2d[M,M]_s^c+\int_0^t
(f(X_s))^2d[A,A]_s^c  +
\sum_{0\leq s\leq t} (F(X_s)-F(X_{s-}))^2.\]

\end{corollary}

{\it Proof}: Follows easily from Theorem 3.2,
 \cite[Theorem, page 144]{F850} and the equality $[X,X]^c=[M,M]^c+[A,A]^c$.
\hfill{\car} \\[2mm]

We are able to prove a version of Theorem 3.2 for weak Dirichlet processes also
with infinite energy. However, in this case we have restricted our attention to
processes with a continuous predictable part.

\begin{theorem}\label{C10}  Let $X=M+A$ be a weak Dirichlet process
with continuous predictable part $A$
and $F$ a
$C^1$ real-valued
function with bounded derivative $f$.

Then the process
$(F(X_t)_{t\geq 0})$ is a
weak Dirichlet process
and the decomposition $F(X)= Y+\Gamma $ holds with the
martingale part
\begin{eqnarray*}& &Y_t= F(0)+ \int_0^tf(X_{s-})dM_s\\
&&\qquad+\int_0^t\int_{\R }\Big(F(X_{s-}+x)-F(X_{s-})-xf(X_{s-})\Big )
(\mu -\nu )(ds, dx)
\end{eqnarray*}

\end{theorem}

{\it Proof}: We consider a sequence $(F^p)_{p\in \N }$ of $C^2$ real functions
such that locally on compact sets $\Vert F-F^p\Vert +\Vert f-f^p\Vert
\rightarrow 0$,
 when
$p\rightarrow
\infty $. Let $A^p$ be a sequence of continuous processes with finite
variation such that
\[\sup_{t\leq T}|A^p_t-A_t|\toP0.\]

Using classical  It\^o's formula for the semimartingale $X^p=M+A^p$
we can write
$$F^p(X^p_t)=F^p(0)+Y^p_t+\Gamma ^p_t$$
where
$$Y^p_t=\int_0^tf^p(X^p_{s-})dM_s+L^p_t$$
with
$$L^p_t=\int_0^t\int_{\R -\{ 0\}
}\Big(F^p(X^p_{s-}+x)-F^p(X^p_{s-})-xf^p(X^p_{s-})\Big )
(\mu -\nu )(ds, dx).$$
Similarly to the proof of Theorem 3.2, we check that

\[\sup_{t\leq T}|Y^p-Y_t|\toP0.\]

On the other hand it is clear that

\[\sup_{t\leq T}|F^p(X^p_t)-F(X_t)|\toP0,\]

which implies that $\Gamma$ as a uniform limit of predictable processes is
also predictable. Finally, by the same arguments as in the proof of Theorem
3.2 we prove that  $[\Gamma,N]=0$
for every continuous local martingale $N$.

\bigskip

\bigskip

\section{Weak  Dirichlet processes and generalized martingale convolutions}

In this section  we deal with processes $X$ such that
\begin{equation}\label{is}
X_t=\int_0^tG(t,s)dL_s
\end{equation}
where $L$ is a quasileft continuous square integrable martingale, and $G$ a
real valued deterministic function of $(s,t)$.

Let us consider the following hypotheses on $G$.

\medskip

\noindent
($H_0$):    $(t,s)\rightarrow G(t,s)$ is
continuous on $\{ (s,t):0<s\leq t\leq T\} $.

\medskip

\noindent
($H_1$):    For all $s$, $t\rightarrow G(t,s)$ has a bounded energy on
$]s,T]$ that is
$$V^2_2(G)((s,T],s)=\sup_n\sum_{t_i\in D_n,t_i\geq s}(G(t_{i+1},s)-G(t_i,s))^2
<\infty .$$

\medskip

\noindent
$$E[\int_0^TV_2^2(G)(]s,T])d[L,L]_s]<\infty \leqno{(H_2):}$$

\medskip

\noindent
$$E[\int_0^T\Gamma ^2(s)d[L,L]_s]<\infty \leqno{(H_3):}$$
where $\Gamma ^2(s)=\sup_{t\leq T}G^2(t,s).$

\medskip

\noindent
\begin{remark}{\rm
Errami and Russo (\cite{ER1}) use, instead of ($H_0$) a slightly
more restrictive assumption, namely: ($H_{0^+}$):    $(t,s)\rightarrow
G(t,s)$ is
continuous on $\{ (s,t):0\leq s\leq t\leq T \} $.

Note that $(H_{0^+})$,  implies  $(H_3)$. Actually
$\Gamma ^2$ is continuous and bounded.

If $t\rightarrow G(t,s)$ admits a quadratic variation on $(s,T]$ along $(D_n)$,
then $(H_1)$ is satisfied.

\medskip

We shall extend $G$ to the square $[0,T]^2$ by setting $G(t,s)=0$ if $s>t$.
}\end{remark}

\begin{theorem} If $X$ meets (\ref{is}) and if $G$ satisfies
$(H_0),(H_1),(H_2),
(H_3)$,

then\begin{description}

\item[{\bf (i)}]
$X$ is a continuous in probability  process with finite energy and
has an optional modification,

\item[{\bf (ii)}] Let us assume that $X$ has a.s.  c\`adl\`ag trajectories,
then $X$ is a weak
 Dirichlet process with natural decomposition $X=M+A$, such that if
$M^n$ is defined as in (\ref{Mn}), then
for every $t\leq T$,

$\vert M^n_t-M_t\vert \rightarrow 0$ in ${\bf L}^2$.

\end{description}
\end{theorem}

\medskip

{\it Proof}: The proof will be given in several steps.

\begin{lemma} $X$ is a continuous in probability process with finite
energy.
\end{lemma}

\medskip

{\it Proof of Lemma 4.1}: First of all, from $(H_2)$ and $(H_3)$, for every
$t\leq T$
$X_t$ is an $\F_t $-measurable square integrable random variable.

  Let us write
$$X_{t_{i+1}}-X_{t_i}=\int_0^{t_i}(G(T_{i+1},s)-G(t_i,s))dL_s+\int_{t_i}^{t_{i+1
}}
G(t_{i+1},s)dL_s.$$
Since $L$ is a square integrable martingale, we get:
\begin{eqnarray*}
&& E(\sum_{t_i\in D_n}(X_{t_{i+1}}-X_{t_i})^2)\cr
&&\leq
2E(\sum_{t_i\in D_n}
(\int_0^{t_i}(G(t_{i+1},s)-G(t_i,s))dL_s)^2)+2E(\sum_{t_i\in D_n}
(\int_{t_i}^{t_{i+1}}
G(t_{i+1},s)dL_s)^2)\cr
&&\leq 2E(\sum_{t_i\in D_n}
\int_0^{t_i}(G(t_{i+1},s-G(t_i,s))^2d[L,L]_s)+2E(\sum_{t_i\in D_n}\int_{t_i}
^{t_{i+1}}(G(t_{i+1},s))^2d[L,L]_s\cr
&&=2E(I^n_1)+2E(I^n_2).
\end{eqnarray*}
By simple calculations
\begin{eqnarray*}I^n_1&=&\sum_{t_i\in D_n}
\sum_{k=1}^i\int_{t_{k-1}}^{t_k}(G(t_{i+1},s)-
G(t_i,s))^2d[L,L]_s\cr
&=&\sum_{t_k\in D_n}
\int_{t_{k-1}}^{t_k}\sum_{i>k}(G(t_{i+1},s)-G(t_i,s))^2d[L,L]_s\cr
&\leq& \sum_{t_k\in D_n}\int_{t_{k-1}}^{t_k}V^2_2(G)((s,T],s)d[L,L]_s,
\end{eqnarray*}
and
\begin{eqnarray*} I^n_2&\leq& \sum_{t_i\in D_n}
\int_{t_i}^{t_{i+1}}(\Gamma ^2(s)d[L,L]_s\cr
&\leq& \int_0^T\Gamma ^2(s)d[L,L]_s. \end{eqnarray*}
Therefore
\begin{eqnarray*} \sup_{n}E(\sum_{t_i\in D_n}(X_{t_{i+1}}-X_{t_i})^2
&\leq& 2E\int_0^TV_2^2(G)((s,T],s)d[L,L]_s+2E\int_0^T(\Gamma (s))^2d[L,L]_s\cr
&<&\infty .\end{eqnarray*}
This proves that $X$ has a finite energy.\\[2mm]

Now, let us take $s,t$ such that $0\leq s<t\leq T$. We get:

$$E[(X_t-X_s)^2]\leq 2E\int_0^s(G(t,u)-G(s,u))^2d[L,L]_u+2E\int_s^t(G(t,u))^2
d[L,L]_u.$$

Since $L$ is continuous in probability, so is $[L,L]$. Under $(H_3)$,
$$2E\displaystyle\int_s^t(G(t,u))^2d[L,L]_u\rightarrow 0,\quad as\,\,
t\rightarrow s,$$
 then by continuity of $t\rightarrow G(t,s)$ and dominated
convergence,
$$E\displaystyle\int_0^s(G(t,u)-G(s,u))^2d[L,L]_u \rightarrow 0.$$
The continuity in probability of the process $X$ follows.

At last, since the process $X$ is $\F_t $-adapted and continuous in
probability, it
admits an
optional modification that we shall denote again by $X$: see for example
\cite{M} pp
230--231,
where Th\'eor\`eme 5 bis is given for a progressively measurable
modification, but the sequence
of approximating
processes introduced in the proof is c\`adl\`ag hence optional (see also below
the proof of the existence of a predictable modification of the process $A$).

Therefore Lemma 4.1 is proven. \hfill{\car}

\begin{lemma} Let us consider the decomposition:

$X_t= A^n_t+M^n_t$,

where as in (\ref{Mn}),
$$M^n_t=  \sum_{t^n_i\in D_n,t^n_{i}\leq t}\left[X_{t^n_{i}}-
E[X_{t^n_{i}}/\F_{t^n_{i-1}}]\right].$$
Then $X$ admits a modification with a decomposition $X_t=A_t+M_t$,
where $M$ is the square integrable martingale
$M_t=\displaystyle\int_0^tG(s,s) dL_s$,
and $A$ a predictable process, such that:

\begin{description}
\item[{\bf(i)}] for every $t\leq T$,

 $\vert M^n_t-M_t\vert \rightarrow 0$ in ${\bf L}^2$,

\item[{\bf(ii)}] for every $t\leq T$, $\vert A_t^n-A_t
\vert \rightarrow 0$ in ${\bf L}^2$,

\item[{\bf(iii)}] for every continuous local martingale $N$, $[A,N]=0$.
\end{description}
\end{lemma}

\noindent

{\it Proof of Lemma 4.2}: For every $t_i, t_{i+1}$ we have
$$E[X_{t_{i+1}}-X_{t_i}\vert
\F_{t_i}]=\int_0^{t_i}(G(t_{i+1},s)-G(t_i,s))dL_s$$
hence for $t\in[0,T]$
$$M^n_t=\sum_{t_{i+1}\leq
t}\int_{t_i}^{t_{i+1}}G(t_{i+1},s)dL_s=\int_0^{\rho ^n(t)}
G^n(s)dL_s,$$
where $G^n(s)=G(t_{i+1},s)$ for $s\in (t_i,t_{i+1}]$ and $\rho ^n(t)=
\max\{ t_{i}: t_{i}\leq t\} $. Note that $\rho ^n(t)\rightarrow t$ as
$n\rightarrow \infty $.

Define $M_t=\displaystyle\int_0^tG(s,s)dL_s$. By $(H_2)$ and $(H_3)$, for
every $\varepsilon >0$
$$E[(M_{\varepsilon })^2]=E\int_0^{\varepsilon }G^2(s,s)d[L,L]_s
\leq E\int_0^{\varepsilon }\Gamma ^2(s)d[L,L]_s$$
and this last expression tends to $0$ when $\varepsilon $ tends to $0$.

Similarly, when $\varepsilon \rightarrow 0$
$$\sup_nE[(M^n_{\varepsilon })^2=\sup_nE\int_0^{\varepsilon }(G^n(s))^2d[L,L]_s
\rightarrow 0.$$

Now, let us fix $\varepsilon >0$ and belonging to $D_n$. Since
$(t,s)\rightarrow G(t,s)$
on $\{ (s,t):\varepsilon \leq s<t\leq T\} $ is continuous, then
it is uniformly continuous and therefore:
$$\sup_{\varepsilon \leq s\leq T}\vert G^n(s)-G(s,s)\vert \rightarrow 0.$$
Hence
\begin{eqnarray*}
E[(M^n_t-M^n_{\varepsilon }-M_t+M_{\varepsilon })^2]&\leq&
2E(\int_{\varepsilon }
^{\rho ^n(t)}(G^n(s)-G(s,s))dL_s)^2\cr
&&+2E(\int_{\rho ^n(t)}^tG(s,s)dL_s)^2\cr
&\leq& 8E\int_{\varepsilon }^1(G^n(s)-G(s,s))^2d[L,L]_s\cr
&&+2E\int_{\rho ^n(t)}^t
G^2(s,s)d[L,L]_s.
\end{eqnarray*}
The first term tends to 0 when $n\rightarrow \infty $
 because $G^n$ converges uniformly to $G$ on $[\varepsilon,T]$.
Now, the continuity in probability of $L$
implies the continuity in probability of $M$ on $[0,T]$, hence the second term
tends to $0$ when $n\rightarrow \infty $.

Note that  for every $n$
$$A^n_t=  \sum_{t^n_i\in D_n, t^n_{i}
\leq t}\E[X_{t^n_{i}}-X_{t^n_{i-1}}/\F_{t^n_{i-1}}]$$
and  $A=X-M$. $A^n$ is a predictable process
and we have $X_{\rho ^n(t)}=A^n_{\rho ^n(t)}+M^n_{\rho ^n(t)}.$

As $X$ is continuous in probability, for every $t\leq T$,
$X_{\rho ^n(t)}\rightarrow X_t$ in ${\bf L}^2$;
moreover, $M^n_t=M^n_{\rho ^n(t)}$ and $A^n_t=A^n_{\rho ^n(t)}$:
we deduce that $A_t^n\rightarrow A_t$ in ${\bf L}^2$ for every $t$.

It follows that $A$ is also adapted and continuous in probability.

Our decomposition $X=M+A$ coincides with the Graversen-Rao decomposition
of Theorem 2.1. But in Theorem 2.1, it is assumed that $X$ is c\`adl\`ag;  here
it is not the case, it is necessary to check that  $A$ admits a
predictable modification.

Actually, since $A$ is continuous in probability on the interval $[0,1]$,
one can find
a subsequence $\{n(k)\}_{k\geq 1}$ such that for every $t\in [0,1]$, $\bar
A^{n(k)}
\rightarrow A_t$ a.s. when $k\rightarrow \infty $, whith $\bar A^n_t=
\sum _i1_{(t^n_i,t^n_{i+1}]}(t)A_{t^n_i}$ and $A^n_0=A_0$. Since every $\bar
A^n$
is a step process adapted and left continuous, it is predictable, and the
process
$A'$ defined by $A'_t=\limsup_k\bar A^{n(k)}_t$ is also predictable. But for
every $t$,
$A'_t=A_t$ a.s. So, we shall suppose now that $A=A'$, and $X=A'+M$.

\medskip

\noindent
{\it Proof of (iii)}: Let $N$ be a continuous local martingale. Using
localization
arguments, we will assume that $N$ is a square integrable martingale. For
every $t$ we can write:
\begin{eqnarray*}
&&\sum_{t_i\in D_n,t_i\leq t}
(N_{t_{i+1}}-N_{t_i})(A_{t_{i+1}}-A_{t_i})\\
& &\qquad\qquad=
\sum_{t_i\in D_n,t_i\leq t}(N_{t_{i+1}}-N_{t_i})
\int_0^{t_i}(G(t_{i+1},s)-G(t_i,s))dL_s\\
&&\qquad\qquad\quad+\sum_{t_i\in D_n,t_i\leq t}
(N_{t_{i+1}}-N_{t_i})\int_{t_i}^{t_{i+1}}(G(t_{i+1},s)-G(s,s))dL_s\\
&&\qquad\qquad=I^n_1+I^n_2
\end{eqnarray*}
Since $N$ is a martingale, using the B-D-G inequality and Schwarz's inequality
we get:
\begin{eqnarray*} E\vert I^n_1
\vert &\leq &cE\left(\sum_{t_i\in D_n}
(N_{t_{i+1}}-N_{t_i})^2
\Bigl(\int_0^{t_i}(G(t_{i+1},s)-G(t_i,s))dL_s\Bigr)^2\right)^{1/2}\cr
&\leq &cE\left(\max_{t_i\in D_n}\vert N_{t_{i+1}}-N_{t_i}\vert
\sum_{t_i\in D_n}\Bigl
(\int_0^{t_i}(G(t_{i+1},s)-G(t_i,s))dL_s\Bigr)^2\right)^{1/2}\cr
&\leq &c(E(\max_{t_i\in D_n}
\vert N_{t_{i+1}}-N_{t_i}\vert ^2))^{1/2}
\left(E\sum_{t_i\in D_n}
\Bigl(\int_0^{t_i}(G(t_{i+1},s)-G(t_i,s))dL_s\Bigr)^2\right)^{1/2}.
\end{eqnarray*}
Because of the continuity of $N$,

$E[\max_{t_i\in D_n}\vert N_{t_{i+1}}-N_{t_i}\vert
^2]\rightarrow 0;$

the second term is estimated as before by

$(E\displaystyle\int_0^TV^2_2(G)((s,T],s)d[L,L]_s)^{1/2}$,

which is finite.

Now, by Schwarz's inequality
\begin{eqnarray*} E\vert I^n_2\vert &\leq& (E\sum_{t_i\in D_n}
(N_{t_{i+1}}-N_{t_i})^2)^{1/2}
\left(E\sum_{t_i\in D_n}
\Bigl(\int_{t_i}^{t_{i+1}}(G(t_{i+1},s)-G(s,s))dL_s\Bigr)^2\right)^{1/2}\cr
&\leq&(E([N,N]_T))^{1/2}\Bigl(E\int_0^T
(G^n(s)-G(s,s))^2d[L,L]_s\Bigr)^{1/2} \end{eqnarray*}
where $G^n$ is defined as above.

Since $E\displaystyle\int_0^T(G^n(s)-G(s,s))^2d[L,L]_s\rightarrow 0$,
we conclude that for every $t$
$$\sum_{t_i\in D_n,t_i\leq t}
(N_{t_{i+1}}-N_{t_i})(A_{t_{i+1}}-A_{t_i})\rightarrow 0$$
in ${\bf L}^1$,
and the covariation process
$[N,A]$ is null for every continuous martingale $N$.

The proof of Theorem 4.1 is complete.\hfill{\car}

\medskip

Unhappily we are not able  to prove that $X$ admits a modification
with c\`adl\`ag trajectories. However, in this direction, we have
 the following lemma.

\begin{lemma}:  $A$  is continuous (hence $X$ is c\`adl\`ag) if the following
additional condition is filled:

\noindent
($H_c$): There exist
$\delta>0, p>1$ and a function $a(u)$ meeting
$$E\Bigl[(\int_0^Ta(u)d[L,L]_u)^p\Bigr]<\infty ,$$
such that for every $s$, $t$, $u$ holds
$$\Bigl(G(t,u)-G(s,u)\Bigr)^2\leq a(u)|t-s|^{{1\over p}+\delta}.$$
\end{lemma}

\medskip

\noindent

{\it Proof }: Let us take $s,t$ such
that $0\leq s\leq t \leq T$. We
 have with a constant $c$ changing from line to line:
\begin{eqnarray*} E(A_t-A_s)^{2p}&\leq& E(\int _0^t(G(t,u)-G(u,u))dL_u-
\int_0^s(G(s,u)-G(u,u))dL_u)^{2p}\cr
&\leq&
cE(\int_0^s(G(t,u)-G(s,u))^2d[L,L]_u)^p\cr
& &\qquad+E(\int_s^t(G(t,u)-G(u,u)^2d[L,L]_u)^p\cr
&\leq& c(t-s)^{1+p\delta }E(\int_0^ta(u)d[L,L]_u)^p\cr
&\leq& c(t-s)^{1+p\delta }.
\end{eqnarray*}
Hence we get the continuity of $A$ by  Kolmogorov's Lemma.\hfill{\car}

\medskip
An analogous result under Holder condition was already given in the paper
(\cite{BM} ), Lemmas 2C and 2D.

\medskip


We shall suppose by now that the processes given by  (\ref{is} ) have a. s.
c\`adl\`ag
trajectories.

\bigskip

We are now interested in investigating conditions on $G$ in order to make
$X$ a Dirichlet process or a weak Dirichlet process admitting a quadratic
variation.  For that let us
consider the following hypotheses:

\medskip

\noindent
$(H_4)$: For all $s$, $t\rightarrow G(t,s)$ has a bounded variation on
$(s,\tau]$,
for every $\tau\leq T$.

(We denote this variation $\vert G\vert ((s,\tau],s)<\infty $).

\noindent
$$ E\int_0^T\vert G\vert ((s,T],s)d[L,L]_s<\infty \leqno{(H_5):} $$

\noindent
$(H_6)$: For all $u, v$,  $t\rightarrow G(t,u)$ and
$t\rightarrow G(t,v)$ have a
finite mutual quadratic covariation with the property $(S)$
on $(\max(u,v),T]$.

(We denote this covariation $ [G(.,u),G(.,v)]_\tau$). Moreover we suppose
that the convergence involved to define the covariation, is uniform
in $u,v$.

\bigskip

Of course $(H_6)$ implies that $t\rightarrow G(t,s)$ admits a quadratic
variation on $(s,T]$ for all $s$ with the property $(S)$ and that $(H_1)$
is satisfied.

\begin{theorem}\begin{description}

\item[{\bf (i)}] Let us assume $(H_0),(H_3),(H_4),(H_5)$. Then $X$ is
a Dirichlet process. (i.e. $A$ is continuous and $[A,A]\equiv 0$).

\item[{\bf (ii)}] Let us assume $(H_0),(H_2),(H_3),(H_6)$. Then $X$ is a
weak Dirichlet
process. Moreover if we assume that the process B defined by
\begin{equation}\label{[A]}
B_t=\int_0^t[G(.,s),G(.,s)]_td[L,L]_s+2\int_0^t(\int_0^v[G(.,u),
G(.,v)]_tdL_u)dL_v \end{equation}
is a c\`adl\`ag process, then $X$ and $A$  admit a quadratic variation such
that:
$[A,A]_t=B_t$
and
\begin{equation} [X,X]=[A,A]+[M,M].
\end{equation}
\end{description}
\end{theorem}

\noindent

\begin{remark}{\rm\begin{description}

\item[{\bf(i)}] If $t\rightarrow G(t,s)$ is $C^1$ for every $s$, denoting
$G_1(t,s)$ its derivative and assuming that $(t,s)\rightarrow G_1(t,s)$ is
continuous on $[0,1]^2$,
we get $A_t=\displaystyle\int_0^t(\int_0^uG_1(u,s)dL_s)du$
by applying Fubini's theorem for stochastic
integrals, so X is a semimartingale.
This result is due to Protter (\cite{P}).

\item[{\bf(ii)}] In case of continuous martingale $L$, part 2) is due to
Errami and Russo
(\cite{ER1} ) and (\cite{ER2} ).
\end{description}
}\end{remark}

\medskip

\noindent

{\it Proof }: (i) First of all, we notice that our hypotheses imply that
$(H_1)$ and $(H_2)$
are satisfied for any sequence $(D_n)$ of subdivisions with mesh
tending to $0$. Since we can write
$$A_{t_{i+1}}-A_{t_i}=\int_0^{t_i}(G(t_{i+1},s)-G(t_i,s))dL_s+
\int_{t_i}^{t_{i+1}}(G(t_{i+1},s)-G(s,s))dL_s,$$
for every $\varepsilon >0$ we get:
\begin{eqnarray*} && E\sum_{t_i\in D_n} (A_{t_{i+1}}-A_{t_i})^2 \\
&& \qquad\leq
E\sum_{t_i\in D_n}\int_0^{t_i}(G(t_{i+1},s)-G(t_i,s))^2d[L,L]_s
\\& &\qquad\qquad+\sum_{t_i\in D_n}\int_{t_i}^{t_{i+1}}
(G(t_{i+1},s)-G(s,s))^2d[L,L]_s\\
&&\leq
\max_{{t_i\in D_n},\varepsilon \leq s\leq T}\vert G(t_{i+1},s)-G(t_i,s)\vert
E\sum_{t_i\in D_n}
\int_{\varepsilon }^{t_i}\vert G(t_{i+1},s)-G(t_i,s)\vert d[L,L]_s\\
&&\qquad\quad
+E\sum_{t_i\in D_n}\int_0^{t_i\wedge \varepsilon
}(G(t_{i+1},s)-G(t_i,s))^2d[L,L]_s \\
&&\qquad\quad+E\int_0^T(G^n(s)-G(s,s))^2d[L,L]_s\\
&&\qquad\leq \max_{{t_i\in D_n},\varepsilon \leq s\leq T}
\vert G(t_{i+1},s)-G(t_i,s)\vert
E\int_{\varepsilon }^1\vert G\vert ((s,1],s)d[L,L]_s\cr
&&\qquad\quad+E\int_0^{\varepsilon
}V^2_2(G)((s,T],s)d[L,L]_s+E\int_0^T(G^n(s)-G(s,s))^2d[L,L]_s
\end{eqnarray*}
Using the following properties:
\begin{eqnarray*}
& &\max_{{t_i\in D_n},\varepsilon \leq s\leq T}
\vert G(t_{i+1},s)-G(t_i,s)\vert
\rightarrow 0,\quad\mbox{\rm
when}\,\,n\rightarrow \infty,\\
& &E\displaystyle\int_{\varepsilon }^T
\vert G\vert ((s,T],s)d[L,L]_s<\infty,
\\& &E\displaystyle\int_0^{\varepsilon }V^2_2(G)((s,T],s)
d[L,L]_s\rightarrow 0,\quad\mbox{\rm when}\,\,
\varepsilon \rightarrow,
0\end{eqnarray*}
and

\[E\displaystyle\int_0^T
(G^n(s)-G(s,s))^2d[L,L]_s\rightarrow 0,\quad\mbox{\rm when }\,\,
n\rightarrow \infty,\]

we deduce that $[A,A]\equiv 0$.
By the inequality
$$\max_{t_i\in D_n}
\vert A_{t_{i+1}}-A_{t_i}\vert^2\leq \displaystyle\sum_{t_i\in D_n}
(A_{t_{i+1}}
-A_{t_i})^2,$$
we get that $A$ is continuous.

\bigskip

(ii) Taking into account Proposition 2.1 we have only to prove the
property $S$ for the process defined by the right hand side of formula
(\ref{[A]}).
Let us notice first that formula (\ref{[A]})
is well-defined, as we take as integrant of
$dL_s$ for the last term the predictable
projection of the optional process
$$(\displaystyle\int_0^s[G(.,u),G(.,s)]_tdL_u)_{s\leq t}.$$
See for example Dellacherie -Meyer (\cite{DM}), Chap.VI for details.

Taking into account that

$E\displaystyle\int_0^T(G^n(s)-G(s,s))^2d
[L,L]_s\rightarrow 0$ when $n\rightarrow \infty $,

we have:

\begin{eqnarray*}
& &[A,A]_t=\lim_n\sum_{t_i\in D_n,t_{i}\leq t}(A_{t_{i+1}}-A_{t_i})^2\cr
&&\qquad=\lim_n\sum_{t_i\in D_n,t_{i}
\leq t}(\int_0^{t_i}(G(t_{i+1},s)-G(t_i,s))dL_s)^2\cr
&&=\lim_n2\sum_{t_i\in D_n,t_{i}\leq
t}\int_0^{t_i}(\int_0^s(G(t_{i+1},u)-G(t_i,u))dL_u)_-(G(t_{i+1},s)
-G(t_i,s))dL_s\cr
&&\qquad\qquad+\lim_n\sum_{t_i\in D_n,t_{i}
\leq t}\int_0^{t_i}(G(t_{i+1},s)-G(t_i,s))^2d[L,L]_s\cr
&&\qquad=2\lim_n I^n_1(t)+\lim_nI^n_2(t). \end{eqnarray*}
Then
$$I^n_2(t)=\sum_k\int_{t_{k-1}}^{t_k}\sum_{i>k,t_{i}\leq t}(G(t_{i+1},s)-
G(t_i,s))^2d[L,L]_s$$

This sequence converges to
$\displaystyle\int_0^t[G(.,s)]_td[L,L]_s$ by dominated convergence.

On the other hand, $I^n_1(t)$ can be written
\begin{eqnarray*}I^n_1(t)&=&\sum_k\int_{t_{k-1}}^{t_k}(\int_0^s\sum_{i>k,
t_{i}\leq t}
(G(t_{i+1},u)-G(t_i,u))(G(t_{i+1},s)-G(t_i,s))dL_u)_-dL_s \cr
&=&\int_0^{\rho ^n(t)}(\int_0^s[G^n(.,u),G^n(.,s)]_tdL_u)_-dL_s.
\end{eqnarray*}
For any optional process $Y$, let us write $Y^P$ its predictable
projection. Then noticing that
$$(\int_0^s[G^n(.,u),G^n(.,s)]_tdL_u)_-=
(\int_0^s[G^n(.,u),G^n(.,s)]_tdL_u)^P$$
and from classical properties of predictable projections, we get for $t\in D_n$
\begin{eqnarray*} &&E(I^n_1(t)-\int_0^t
(\int_0^s[G(.,u),G(.,s)]_tdL_u)dL_s)^2 \cr
&&\qquad=
E[\int_0^t[(\int_0^s([G^n(.,u),G^n(.,s)]_t-[G(.,u),G(.,s)]_t)dL_u)^P]^2
d<L,L>_s] \cr
&&\qquad\leq
E[\int_0^t[(\int_0^s([G^n(.,u),G^n(.,s)]_t-[G(.,u),G(.,s)]_t)dL_u)^2]^P
d<L,L>_s] \cr
&&\qquad=E[\int_0^t(\int_0^s([G^n(.,u),G^n(.,s)]_t-[G(.,u),G(.,s)]_t)dL_u)^2
d<L,L>_s] \end{eqnarray*}

We show now that for any $t\in \bigcup_nD_n$

\[\int_0^t(\int_0^s([G^n(.,u),G^n(.,s)]_t-[G(.,u),G(.,s)]_t)
dL_u)^2d<L>_s\toP0,\quad\mbox{ \rm when}\,\,n\rightarrow \infty.
\]

Note that
\begin{eqnarray*}
&&E(\int_0^s([G^n(.,u),G^n(.,s)]_t-[G(.,u),G(.,s)]_t)dL_u)^2 \cr
&&\qquad=E(\int_0^s([G^n(.,u),G^n(.,s)]_t-[G(.,u),G(.,s)]_t)^2d<L,L>_u).
\end{eqnarray*}
But for every $u,s,t$, $[G^n(.,u),G^n(.,s)]_t-[G(.,u),G(.,s)]_t\rightarrow 0$
when $n\rightarrow \infty $, and we have the estimation
\begin{eqnarray*}&&
\vert [G^n(.,u),G^n(.,s)]_t-[G(.,u),G(.,s)]_t\vert \cr
&&\qquad\leq 1/2([G^n(.,u)]_t+[G^n(.,s)]_t+[G(.,u)]_t+G(.,s)]_t).
\end{eqnarray*}
>From ($H_6$) this last term is bounded, hence by dominated convergence
$$E(\int_0^s([G^n(.,u),G^n(.,s)]_t-[G(.,u),G(.,s)]_t)dL_u)^2\rightarrow 0$$
for every $s,t$. So, we can get easily (18) by localisation of $L$.

We are finished as soon as we remark that using the continuity of process
$<L,L>$,
for every $t$
\begin{eqnarray*} &&E(I^n_1(t)-I^n_1(\rho ^n(t))-\int_{\rho ^n(t)}^t
(\int_0^s[G(.,u),G(.,s)]dL_u)dL_s)^2 \cr
&&\qquad\qquad= E(\int_{\rho ^n(t)}^t (\int_0^s[G(.,u),G(.,s)]dL_u)dL_s)^2 \cr
&&\qquad\qquad=E\int_{\rho ^n(t)}^t (\int_0^s[G(.,u),G(.,s)]dL_u)^2d<L,L>_s
\end{eqnarray*}
converges to $0$ when $n\rightarrow \infty $. \hfill{\car}

\bigskip

\noindent
{\bf Example 1: Fractional normal processes of index $H>1/2$}

We consider the case where $L$ is a normal martingale (i.e. a square
integrable martingale
with predictable quadratic variation $<L,L>_t=t$), and $G(t,s)$ given for
$t\geq s$ by:
\begin{equation} G(t,s)=cs^{1/2-H}\int_s^tu^{H-1/2}(u-s)^{H-3/2}du
\end{equation}
with $c$ constant. Of course, when $L$ is a standard Brownian motion, $X$
is the classical fractional Brownian motion.

Let us check that ($H_3$),($H_4$) and ($H_5$) are satisfied.
\begin{eqnarray*} \vert G\vert ((s,T],s)&=&
cs^{1/2-H}\int_s^Tu^{H-1/2}(u-s)^{H-3/2}du \cr
&\leq& cs^{1/2-H}\int_s^T(u-s)^{H-3/2}du \cr
&=&cs^{1/2-H}{1\over {H-1/2}}(1-s)^{H-1/2}du \cr
&\leq& {c\over {H-1/2}}s^{1/2-H} \end{eqnarray*}
and
$$\int_0^T\vert G\vert ((s,T],s)ds\leq {c\over {(H-1/2)(3/2-H)}},$$
hence finally
$$\int_0^T\vert G^2\vert ((s,T],s)ds\leq {c\over {(H-1/2)(2-2H)}}.$$

This means that, for every normal martingale $L$, the process $X$ defined by
$X_t=\displaystyle\int_0^tG(t,s)dL_s$ is a Dirichlet process. Since
$G(s,s)=0$, its
 martingale part is null. In particular $X$ is continuous, even if $L$ is not.

\bigskip

\noindent
{\bf Example 2: Weak Dirichlet process driven by a Brownian motion}

For this example we take $L=B$ a standard Brownian motion and $G(t,s)=\beta
(t)f(s)$ for
$t\geq s$, where $t\rightarrow \beta (t)$
is a fixed Brownian trajectory such that its quadratic variation is
$[\beta (.)]_t=t$, and $f$ is a real continuous function on $[0,T]$.

Here we can apply part 2) of Theorem 4.2. Actually:
$$[G(;,u),G(.,v)]_t= \lim_{t_{i+1}\leq t,t_i\geq max\{ u,v\} }
\sum_i(\beta (t_{i+1})-\beta (t_i))^2f(u)f(v)$$
and this term converges uniformly in $u,v$ to $(t-\max\{ u,v\} )f(u)f(v)$.

Therefore we get the decomposition $X=M+A$, with
\begin{equation} M_t=\int_0^t\beta (s)f(s)dB_s \end{equation}
and the formula (\ref{[A]}) gives the quadratic variation of $A$
\begin{equation} [A,A]_t=\int_0^t(t-s)f^2(s)ds+2\int_0^t\int_0^s(t-s)f(u)f(s)
dB_udB_s \end{equation}
which can be written
$$[A,A]_t=\int_0^t(\int_0^sf(u)dB_u)^2ds.$$
In particular the process
$(\displaystyle\int_0^t\int_0^s(t-s)f(u)f(s)dB_udB_s)_t$ has
a finite variation.

Since $E([A,A]_t)= \displaystyle\int_0^t(t-s)f^2(s)ds\not = 0$, $X$ is not a
Dirichlet process; however, Theorems 4.1 and 4.2 ensure that $X$
is a weak Dirichlet process admitting quadratic variation.

\bigskip
Actually, this example 2 is a particular case of the following:

\bigskip
\noindent
{\bf Example 3}

Let us consider $G(t,s)$ under the form

$G(t,s)= \int_s^tf(u,s)d\beta _u$

\noindent
where for every $s<t$, $u\rightarrow f(u,s)$ has a bounded variation on
$(s,t]$, and
$(u,s)\rightarrow f(u,s)$ is continuous on $\{(u,s): 0\leq u\leq s\leq T\}$. We
denote $df_u(u,s)$ the measure associated to the variation process. We
assume also that
$\beta $ is a deterministic real continuous on $[0,T]$ function admitting a
quadratic
variation along $(D_n)$, and $\beta _0=0$.

 $X_t=A_t$ and we shall prove the Fubini type formula :

\begin{equation} A_t=\int_0^t(\int_s^tf(u,s)d\beta
_u)dL_s=\int_0^t(\int_0^uf(u,s)dL_s)
d\beta _u.
\end{equation}

Indeed, $A$ admits a quadratic variation along $(D_n)$ given by

$$[A,A]_t=\int_0^t(\int_0^sf(u,s)dL_u)^2d[\beta ,\beta ]_s.$$
\medskip

Actually, taking into account that $[\beta ,f(.,s)]=0$, we get

$$\int_s^tf(u,s)d\beta _u=\beta _tf(t,s)-\beta _sf(s,s)-\int_s^t\beta _u
df_u(u,s).$$

Then

$$A_t= \beta _t\int_0^tf(t,s)dL_s-\int_0^t\beta _sf(s,s)dL_s-
\int_0^t(\int_s^t\beta _udf_u(u,s))dL_s.$$

>From Theorem 4.2 (i) the process $Y$ defined by $Y_t=
\displaystyle\int_0^tf(t,s)dL_s$ is a
Dirichlet process and $[\beta ,Y]=0$, then by integration by parts, we get;

$$\beta _tY_t= \int_0^tY_sd\beta _s+\int_0^t\beta _sdY_s$$

And by using the  sequence $(D_n)$

\begin{eqnarray*}&\int_0^t\beta _sdYs= lim_{D_n}\sum_{t_i \in D_n, t_i\leq
t}\beta _{t_i}
(Y_{t_{i+1}}-Y_{t_i})\cr
&=\lim_{D_n}\sum_{t_i \in D_n, t_i\leq t}\beta
_{t_i}\int_{t_i}^{t_{i+1}}f(t_{i+1},s)dL_s\cr
&+\lim_{D_n}\sum_{t_i \in D_n, t_i\leq t}\beta
_{t_i}\int_0^{t_i}(f(t_{i+1},s)-f(t_i,s))dL_s \cr
&=\int_0^t\beta _sf(s,s)dL_s+\int_0^t(\int_s^t\beta _udf_u(u,s))dL_s,
\end{eqnarray*}
then we deduce the formula (18).

\section{Appendix}
\subsection{A process with finite energy but without quadratic variation
along the
dyadics.}

Let $D_k$ be the $k$-th dyadic subdivision of $[0,1]$, that is
$\displaystyle t^k_j={j\over 2^k}$,$0
\leq j\leq 2^k$.

We are willing to build a deterministic function $x$ such that $S_k=1$ if
$k$ is even and greater
 than 2, and $S_k=2$ for $k$ odd. Such a funcion has obviously finite
energy along the
sequence $(D_k)_k$ (and indeed its energy is equal to 2, although is would
be equal to 1 along the
sequence $(D_{2k})_k$), but has no quadratic variation since the sequence
$(S_k)$ has 2 accumulation points.

Let us begin with defining $x_0=x_1=0$ and $x_{1/2}=1$, so that $S_1=2$.

At the second step, we define $x_{1/4}=x_{3/4}=1/2$, so that $S_2=1$

In order to make our construction clear, we go into details for the third step.

We want to define $x_{j/8}$ for odd $j$ in order that $S_3=2$. The idea is
to compute $x_{j/8}$ such
that
$$\Bigl(x_{j \over 8}-x_{{j-1}\over 8}\Bigr)^2+\Bigl(x_{{j+1} \over
8}-x_{{j}\over 8}\Bigr)^2=2\times
 \Bigl(x_{{j+1} \over 8}-x_{{j-1}\over 8}\Bigr)^2.$$
Actually, this amounts to find a solution $y$ to an equation like
\begin{equation}\label{eq2}
(a-y)^2+(b-y)^2=2\times (a-b)^2.
\end{equation}

As equation (\ref{eq2}) has two solutions, namely $((1+\sqrt 3)a+(1- \sqrt
3)b)/2$ and $((1-\sqrt 3)a+(1+ \sqrt
3)b)/2$, we have 2 possible choices for each $x_{j/8}$ with odd $j$ in order
that $S_3=2$.

This process is then iterated as follows :

1. Assume that we have constructed $x_{j\over 2^{2k-1}}$, $0\leq j\leq
2^{2k-1}$,
such that  $S_{2k-1}=2$ for some $k$. Then we put $x_{2j+1\over
2^{2k}}=(x_{2j\over 2^{2k}}+x_{2j+2\over 2^{2k}})/2$ (so that it is the middle
of its neighbours). Then it is readily checked that $S_{2k}=1$.

2. Now we have to choose the $x_{2j+1\over
2^{2k+1}}$'s. We will proceed as was done above for $k=1$. Namely, we can
always
choose $y=x_{2j+1\over
2^{2k+1}}$ so that it solves equation (\ref{eq2}) with $a=x_{2j\over
2^{2k+1}}$ and $b=x_{2j+2\over
2^{2k+1}}$, and the result follows the same lines as for $k=1$.

\medskip

It remains to check that we can build a real continuous function $x$ on $[0,1]$
with the specified values on the dyadics.

Let $x^n$ be the piecewise linear function joining the points
constructed at rank $n$. We will show that the sequence $(x^n)$ satisfies a
uniform Cauchy criterion,
which will give the claim.

\medskip

First note that it is obvious (again from the
solution of equation (\ref{eq2}) that any two
neighbours at rank $2k$ or $2k+1$ are far from each other at most $((1+\sqrt
3/4)^k$. In other words, we have always
\begin{equation}\label{voisins}
|x^n_{i+1\over2^n}-x^n_{i\over 2^n}|\leq \Bigl({1+\sqrt 3\over
4}\Bigr)^{n\over 2}.
\end{equation}

  Now, fix $\varepsilon>0$. For positive $n$ and $p$ and for $t\in[0,1]$,
let $t^n_i$
be the closest to $t$ point in $D_n$, and $t^{n+p}_j$ the closest to $t$
point in $D_{n+p}$.
  Without loss of generality we will assume that $t^n_i\leq t^{n+p}_j \leq
t^n_{i+1}$.
  We have then

  \begin{eqnarray}\label{ineq}
  |x^n_t-x^{n+p}_t|&\leq
|x^n_t-x^{n}_{t^n_i}|+|x^{n}_{t_i^n}-x^{n+p}_{t_i^n}|\cr
  &+|x^{n+p}_{t_i^n}-x^{n+p}_{t^{n+p}_j}|+|x^{n+p}_{t^{n+p}_j}-x^{n+p}_t|.
 \end{eqnarray}

>From (\ref{voisins}) and the definition of $x^n$ it is obvious that the
first and last terms in the
right-hand side of (\ref{ineq}) can be made as small as wanted (say less
than $\varepsilon/3$),
uniformly in
$t$ and $p$, for $n$ large enough. Moreover, as $D_n\subset D_{n+p}$, the
second term is identically zero.

Hence it remains to uniformly estimate
$|x^{n+p}_{t_i^n}-x^{n+p}_{t^{n+p}_j}|$ (note that in principle,
 at most $2^p$ points of $D_{n+p}$ may lay between $t_i^n$ and $t^{n+p}_j$,
so that there is no trivial uniform in $p$ majorization of such a sequence).

Let us put $t_n^i=k/2^n$

We assume for instance that $l:=x_{k/2^n}\leq h:x_{(k+1)/2^n}$ and that $n$
is odd. Then clearly $l
\leq x_{(2k+1)/2^{n+1}}\leq h$ (since
$x_{(2k+1)/2^{n+1}}=(l+h)/2$), and if we choose
$x_{(4k+1)/2^{n+2}}=((1+\sqrt 3)l+(1- \sqrt
3)(l+h)/2)/2$ and $x_{(4k+3)/2^{n+2}}=((1+\sqrt 3)(h+l)/2+(1- \sqrt
3)h)/2$, we can see that both values lay in $[l-(\sqrt 3 -1)(h-l)/4, h]$.
Keeping
(\ref{voisins}) in mind,
 we conclude that for every point $s$ in $D_{n+2}\cap
[k/2^n, {(k+1)/ 2^n}]$, and for every $p \geq 2$,

$$l-(\sqrt 3 -1)/4\Bigl({1+\sqrt 3\over 4}\Bigr)^{n\over
2}\leq x^{n+p}_s\leq h.$$

If we iterate this procedure, it is straightforward now that for every
point $s$ in $D_{n+2m}\cap
[k/2^n, {(k+1)/ 2^n}]$, and for every $p \geq 2m$,

$$l-{{\sqrt 3-1}\over 4}\Bigl({1+\sqrt 3\over 4}\Bigr)^{n\over
2}\sum_{i=0}^{m-1}
\Bigl({1+\sqrt 3\over 4}\Bigr)^i\leq x^{n+p}_s\leq h $$

and eventually, for every point $s$ in $\bigcup_m D_{n+2m}\cap
[k/2^n, {(k+1)/ 2^n}]$, and for every $p \geq 0$,

$$l-{{\sqrt 3 -1}\over 4}\Bigl({1+\sqrt 3\over 4}\Bigr)^{n\over 2}{4\over
3-\sqrt
3}\leq x^{n+p}_s\leq h.
$$

At last, il follows that the third term in the right-hand of (\ref{ineq})
can be made as small as wanted,
 say less than $\varepsilon/3$, for $n$ big
enough,uniformly in $t$ and $p$, and finally we checked the uniform Cauchy
criterion for the sequence of
 functions $(x^n)$. Hence this sequence converges to a
continuous function $x$, such that by construction, for every $k$

$$\sum_{i\in D_{2k}}(x_{i+1}-x_i)^2=1$$

and

$$\sum_{i\in D_{2k+1}}(x_{i+1}-x_i)^2=2.$$

This function has finite energy, but no quadratic variation along the dyadics.

\subsection{A continuous function with discontinuous pre-quadratic variation}
We consider the function introduced in \cite{CS}, Example 1, that is the
piecewise affine function $X$
such that $X_t=0$ at each $t=1-2^{1-2p}$, $X_t=1/p^{1/2}$ for
$t=1-2^{-2p}$, and $X$ is affine between these points. If we define
moreover $X_1=0$, $X$ is a
continuous fonction on $[0,1]$.

It is clear that $X$ is a fonction of finite variation, hence of zero
quadratic variation on every
$[0,t]$ with $t<1$.

On the other hand, it was proven in \cite{CS} that $X$ has an infinite
quadratic variation on $[0,1]$ along $S:=\{1-2^{-2k}, k\geq 1\}$.

For $n>0$, we define now a subdivision $\pi_n$ of $[0,1]$ as follows :
$$\pi_n=\bigcup_{j\leq 2^{2n}-1} \{j2^{-2n}\}\cup\bigcup_{k\geq
n}\{1-2^{-2k}\}.$$
It is straightforward that for every $n$ $X$ has an infinite quadratic
variation along $\pi_n$,
although its quadratic variation along $\pi_n\cap[0,t]$ goes to zero as
$n\to\infty$ for
every $t<1$.

\medskip
Note that if we modify our example in order that $X_t=1/p$ for
$t=1-2^{-2p}$, $p>0$, everything else remaining unchanged, the
pre-quadratic variation of $X$ on
$[0,t]$ is equal to zero if $t<1$, but finite and non zero for $t=1$.

\subsection{A continuous in probability process may admit no c\`adl\`ag
modification.}

Let $X$ be the piecewise affine function such that $X_t=0$ at each
$t=1-2^{1-2p}$, $X_t=1$ for
$t=1-2^{-2p}$, and $X$ is affine between these points. If we define
moreover $X_t=0$ outside $[0,1)$,
 we get a discontinuity of the second kind at 1.

Now define the process $Y$ as follows : $Y_t=X_{t-T}{\bf 1}_{t\geq T}$,
where $T$ is a random variable
 uniformly distributed on $[0,1]$, then $Y$ is continuous in probability,
but every path of $Y$ has
almost surely a discontinuity of the second kind between times 1 and 2.

\medskip
Note that this result remains true even if we ask our process to have a
quadratic variation along a
sequence $(\pi_n)$.

\end{document}